\documentclass[12pt, a4paper]{article}
\usepackage[utf8]{inputenc}
\usepackage{graphicx}
\usepackage{amssymb}
\usepackage{stmaryrd}
\usepackage{enumitem}

\usepackage{geometry}
\usepackage{tikz}
\usepackage{float}
\usepackage{url}
\usepackage{mathtools}
\usepackage{amsmath,amsfonts,amsthm}
\usepackage{mathrsfs}

\usepackage{amsthm}

\theoremstyle{plain}

\theoremstyle{definition}

\theoremstyle{remark}

\usepackage{pgfplots}
\pgfplotsset{compat=1.18}

\usepackage{subfig}

\usepackage{tikz}
\usetikzlibrary{matrix}
\usepackage{xcolor}

\usepackage{ytableau}
\usepackage[colorlinks,
linkcolor=blue,
anchorcolor=blue,
citecolor=blue
]{hyperref}
\begin{document}
\begin{center}
{\large \bf The Orientalists' Stance Towards Arabic Sciences (Especially Arabic Astronomy)}
\end{center}
\begin{center}
 Duaa Abdullah$^1$ \hspace{0.2 cm} \& \hspace{0.2 cm} Jasem Hamoud$^{2}$\\[6pt]
 $^{1,2}$ Physics and Technology School of Applied Mathematics and Informatics \\
Moscow Institute of Physics and Technology, 141701, Moscow region, Russia\\[6pt]
Email: $^{1}${\tt duaa1992abdullah@gmail.com},	 $^{2}${\tt hamoud.math@gmail.com}
\end{center}
\noindent
\begin{abstract}
In this paper, we highlight the influence of Arab/Islamic civilization in the field of the history of astronomy on European historians. 
We also aim to elucidate the stance of Orientalists toward the study of Arab sciences and to clarify their orientations, with a particular focus on astronomy, while revealing the significant role played by Arab scholars in this domain and the impact of their contributions—especially astronomical tables (zij)—on Western astronomers. 
Furthermore, we have clarified the mechanisms of transmission of Arab sciences, particularly astronomy, from Arab scholars to Western scholars, and the role of Arab astronomers in Western civilization. 
In addition, we address the contributions of Arab scholars to the development of astronomy and the perspective of Orientalists, particularly David King, regarding this matter. We also underscore the importance of Orientalists’ works in analyzing Arab/Islamic scholarly output, identifying its influence on the West in the field of astronomy, and demonstrating how Western scholars benefited from translations of Arabic books in this discipline. In this paper, we adopt the historical retrieval methodology, by referencing previously documented astronomical information and contributions, with an emphasis on the processes of transmission of these sciences from the Arabs to the West.
\end{abstract}

\noindent\textbf{MSC Classification 2010:} 05C05, 05C12, 05C20, 05C25, 05C35, 05C76, 68R10.

\noindent\textbf{Keywords:} Arabic sciences, Astronomy, Islamic civilisation, Arabs, Orientalism, Orientalists, Astronomical tables.

\section{Introduction}\label{sec1}

The topic of Orientalism and its role in the study of Arabic sciences constitutes an intertwined and complex field that reflects the cultural and historical interaction between the East and the West. European historians exerted tremendous efforts in the domain of astronomy, not confining themselves to European astronomy alone; rather, they endeavored to comprehend and trace the evolution of this science across other civilizations through centuries and millennia~\cite{Chaouch2001}.

The influence of Arab/Islamic civilization on the Western world was profoundly significant, as numerous astronomical works—particularly Arabic astronomical tables (zijes)—were transmitted to them. The most compelling evidence lies in the acknowledgment by many Orientalists of the monumental role played by Arab/Islamic culture in enriching European thought. Astronomy owes a substantial portion of its advancement to the Arabs, whose contributions played a pivotal role in its development, preservation, and transmission to the West through translation, refinement, correction via precise observation and verification, and subsequent conveyance.

Orientalism evolved over time to become a research subject necessitating distinguished cultural and linguistic knowledge, coupled with the capacity for critical and rigorous analysis of Arab and Islamic works. Nevertheless, it was not unaffected by political ideological cultures that managed to alter its interpretations and outcomes.

\begin{itemize}
    \item \textbf{Significance of the Research:} To highlight the stance of Orientalists toward Arabic sciences and to elucidate the true extent of the substantial role played by Arab astronomers in contributing to Western astronomy, particularly their contributions to astronomical tables (\textit{Z\={\i}j}es), and their influence and role therein.
    \item \textbf{Objectives of the Research:} To identify the methods of transmission of Arabic sciences, especially astronomy, to the West through Arab scholars.  To underscore the role of Arab astronomers in Western civilization.
    \item \textbf{Research Methodology:} The historical (retrospective) method has been adopted by referencing previously documented astronomical information and \textit{Z\={\i}j}es, citing scholars and explaining how this knowledge was transferred from the Arabs to the West.
    \item \textbf{Definition of Orientalists and the Scope of Their Work:}  The origins of Orientalism trace back to the Middle Ages, when leading European universities began offering orientations to delve deeply into Arab/Islamic civilization due to its significant contributions across numerous fields, such as the sciences, philosophy, and others.
   Orientalists~\cite{Tamimi2021} are Western scholars and researchers who specialized in studying the cultures and civilizations of the Middle East, with particular emphasis on Islamic civilization and Arabic sciences. The term derives from \emph{Orientalism}, which denotes an attraction to Eastern studies.
   
   Orientalists play a crucial role in analyzing Eastern texts, translating them, and understanding their cultural and intellectual components within their historical context. Their fields of study encompass diverse aspects, including language, literature, philosophy, history, and the natural sciences, thereby illuminating Arab and Islamic achievements in human civilization and the history of scientific knowledge. 
\end{itemize}

\section{Influence of Arab Scholars' Works on the West in the Field of Astronomy}
   Astronomy was known in Iraq more than a century before the Arab conquest according to~\cite{Wat2016}. Muslim scholars advocated for abolishing astrology, which was based on illusion, and instead directed the science of the stars toward truths grounded in observation, experimentation, and empirical knowledge. 
   They observed the celestial spheres, compiled \textit{Z\={\i}j}es, measured latitudes, and monitored the planets. 
   They delved into the works of earlier scholars, completed what was lacking, and synthesized them. 
   
   Thus, they achieved considerable success in astronomy~\cite{Farraj2002} in particular and made significant advancements to it. 
   Arab scholars acquired astronomical knowledge from the Indians, Persians, and Greeks, refined it, corrected it, expanded its astronomical discussions~\cite{Yagi1997}, presented it in a new form, and subjected every issue to rigorous investigation and scrutiny~\cite{Badr1982}. 
   
   They authored numerous books and clarified the computational basis for each step, rendering astronomy markedly different from its earlier state.
   When astronomy garnered attention from the Arabs, they began translating Sanskrit\footnote{\textbf{Sanskrit:} an ancient language of India; and Pahlavi: Middle Persian, which evolved across multiple eras. Both are regarded as pivotal sources upon which Arab astronomers relied in their scholarly works.}, Pahlavi, Greek, and Syriac works~\cite{Wat2016}.
   
   Among the most important achievements of Arabs and Muslims in astronomy were the compilation of \textit{Z\={\i}j}es and the preservation of ancient heritage from loss by translating astronomical books into Arabic, such as Ptolemy's Almagest (Ptolemy, ca. 100–170 CE; Arabic dates erroneously given as (554 BC/90CE–484 BC/168 CE)~\cite{Mahfouz2010,Hamd2001}, whose original Greek text was lost, surviving only in its Arabic translation~\cite{Mahasna2001}.
   
   The Almagest served as the foundational theoretical text, translated by the late 8th century CE, revised multiple times, and accompanied by numerous commentaries and introductions~\cite{Wat2016}. The first to oversee its interpretation and translation into Arabic was the Barmakid Yahya ibn Khalid (d. 805 CE/190 AH), minister to Harun al-Rashid~\cite{Sezgin2008}, or alternatively through al-Hasan ibn Quraysh (9th century CE/3rd century AH) for Caliph al-Ma'mun~\cite{Nadim1971,King2001}.
   
   It remained a reference until the 10th century CE, when the scholar Abd al-Rahman al-Sufi (903–986 CE/291–376 AH) translated it~\cite{Sufi1981}, added his personal observations, and corrected its errors, producing his Book on the Images of the Forty-Eight Constellations, which was more accurate and comprehensive than the Almagest and served as a reference in astronomy until the early 20th century~\cite{Sawaf85}.
   
   Arab astronomers, following Ptolemy's view, believed in a stationary Earth with eight spheres revolving around it: the Sun, Moon, five planets, and fixed stars. Applying this system to observed phenomena required numerous astronomical instruments. Over time, the Arabs recognized the weaknesses in Ptolemy's theory and criticized it, despite the significant simplifications to astronomy by Ibn al-Shatir al-Dimashqi (ca. 1350 CE/750 AH). Most Arab astronomical works were not merely theoretical but focused on \textit{Z\={\i}j}es of Indian, Persian, and Greek origin. Due to inconsistencies among Arab astronomers' tables, more accurate results were achieved. Al-Battani's tables from ca. 900 CE/286 AH were of utmost precision; his eclipse observations were utilized as late as 1749 CE/1162 AH~\cite{Wat2016}.
   
   In the late Middle Ages and early modern period, Europe drew from Arabic science and taught the books authored by those Arab geniuses. Europeans thus derived the foundations of astronomy from Arabic references, through which the West became acquainted with Ptolemy and his \emph{Almagest}, a title bestowed by the Arabs~\cite{Yagi1997}.
   When the telescope was invented in the early 17th century CE/11th century AH by the eminent scientist Galileo Galilei, new horizons opened for astronomy.
   
   Al-Biruni's \emph{Al-Qanun al-Mas'udi} (973–1048 CE/362–440 AH) became the primary source for all astronomical works~\cite{Saadi2012,Hamd2001}, including the \textit{Z\={\i}j}es of al-Battani (858–929 CE/244–317 AH)~\cite{Zirkali2002v6}, al-Tusi (1201–1274 CE/597–672 AH)~\cite{Zirkali2002v8}, and others. Al-Sufi's Book on the Images of the Forty-Eight Constellations served as the fundamental celestial atlas.
   When the Western astronomer Hevelius observed the Moon—and was the first to do so in 1647—he identified several craters and named some after Arab astronomers. Ten craters bear the names of Arab astronomers, as follows:
\begin{enumerate}
\item Al-Farghani (d. ca. 1007 CE/398 AH)~\cite{Zirkali2002v1}.  
\item Al-Battani.  
\item Al-Sufi; fixed stars were derived from al-Sufi, retaining the names he assigned.  
\item Al-Hasan ibn al-Haytham (965–1038 CE/354–430 AH)~\cite{Zirkali2002v6}. 
\item  Nasir al-Din al-Tusi. 
\item Thabit ibn Qurra (836–901 CE/221–288 AH)~\cite{Zirkali2002v2}.
\item Ulugh Beg (1393–1449 CE/796–853 AH)~\cite{AlMumin1992}.
\item Jabir ibn Aflah (d. ca. 1145 CE/450 AH)~\cite{AlMumin1992}.  
\item Al-Ma'mun (814–833 CE/198–217 AH), in recognition of his patronage of astronomy~\cite{King2001}. 
\item Abu al-Fida' (1273–1331 CE/672–732 AH)~\cite{AbuAlFida}.
\end{enumerate}
These are clearly visible on the full segmented lunar map divided into two inverted sections, as images typically appear inverted in astronomical observatories, as shown in Figure~\ref{fignduaan1} below~\cite{Badr1982}:
\begin{figure}[H]
    \centering
    \includegraphics[width=0.5\linewidth]{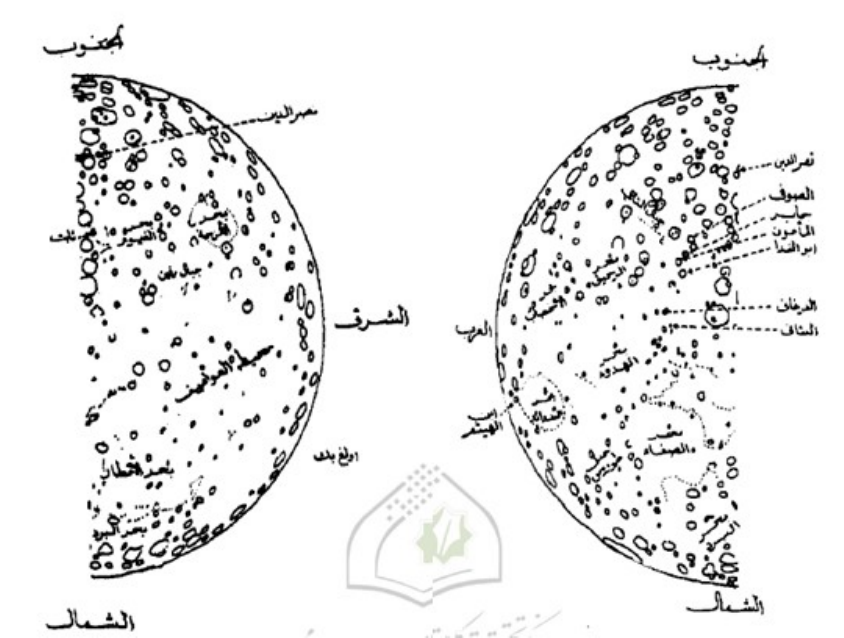}
    \caption{ The complete lunar map divided into two mirrored halves (the eastern half and the western half)~\cite{Badr1982}.}
    \label{fignduaan1}
\end{figure}
A substantial number of star names and astronomical terms found in modern foreign languages are directly derived from ancient Arabic vocabulary~\cite{Farroukh1970}.

The foundational advancement of Arab scholars in the \emph{mathematical sciences} was instrumental in establishing their preeminence in \emph{astronomy}, a field they cultivated with immense dedication. This commitment is demonstrably evidenced by the numerous sophisticated astronomical observatories established across the Islamic world, including those in Samarkand, Damascus, Cairo, Fez, Toledo, Cordoba, and Maragheh. Arab astronomers utilized highly precise instruments for observation, such as the \emph{astrolabe}, the sundial (\emph{Mizwala}), and the altitude measuring instrument, among others.

Among the most distinguished figures in this tradition was \emph{Abu al-Wafa' al-Buzjani} (328 AH/940 CE – 387 AH/997 CE)~\cite{Brockelman59}, who compiled comprehensive astronomical tables in his \emph{Al-\textit{Z\={\i}j} al-Wadih} (The Clear Tables), which meticulously calculated the positions of stars. Equally influential was \emph{Al-Battani}, the author of \emph{Al-\textit{Z\={\i}j} al-Sabi'} (The Sabean Tables). Both of these \emph{\textit{Z\={\i}j}} (astronomical handbooks) had a profound and lasting impact on the development of astronomical science in both the Eastern and Western worlds.

Arab scholars achieved remarkable scientific breakthroughs, including: establishing the \emph{sphericity of the Earth} and its rotation around its axis, accurately estimating its circumference, and precisely measuring \emph{the solar year}. Furthermore, they compiled detailed tables for the positions of the planets, utilizing their advanced astronomical instruments, such as the astrolabe~\cite{Yagi1997}, for their observations~\cite{Hamoud2023}.

Astronomy was one of the first sciences to receive significant patronage in Baghdad. The study of its principles was not confined to Arab scholars; subsequent generations followed their methodology and inherited their knowledge. This lineage is exemplified by \emph{Ulugh Beg }and his famous \emph{\textit{Z\={\i}j}}, which can be considered the final, exemplary product of the Baghdad school, whose golden age spanned seven centuries (132 AH/750 CE – 853 AH/1450 CE).

The transmission of Indian knowledge to the Arab world occurred through the translation of the book \emph{Sindhind} into Arabic during the reign of Caliph Al-Mansur. This translation was initially performed by \emph{Al-Fazari} (who flourished in the 2nd Hijri century/8th century CE)~\cite{AlQiftiND}, and subsequently by \emph{Al-Khwarizmi }(died after 232 AH/847 CE)~\cite{Brockelman59}. Muhammad al-Fazari's translation of \emph{Sindhind} for Al-Mansur served as the foundational text for astronomical calculations among the Arabs and remained their primary reference until the era of Al-Ma'mun.

During Al-Ma'mun's time, \emph{Muhammad ibn Musa al-Khwarizmi} excelled. Possessing extensive knowledge of the stars, he created his own \emph{\textit{Z\={\i}j}}, the \emph{\textit{Z\={\i}j} al-Khwarizmi}, which synthesized the astronomical traditions of the Indians, Persians, and Romans. While he based its foundation on the \emph{Sindhind}, he introduced his own corrections (\emph{ta'adil}) and obliquity (\emph{mayl}). His corrections followed the Persian school, and the obliquity of the sun was based on the Ptolemaic system. The work was structured into excellent chapters, earning the admiration of his contemporaries, although its chronology was based on the Persian calendar.

This work was later adapted to the Arabic calendar by \emph{Maslama ibn Ahmad al-Majriti al-Andalusi} (d. 398 AH/1007 CE)~\cite{Farraj2002}. Through this process, the Arabs became intimately acquainted with the Indian systems of mathematics, studied the Indian numerals, refined them, and developed two distinct series:
\begin{enumerate}
    \item The first, known as the \emph{Hindu numerals} (or Eastern Arabic numerals), comprising: "1, 2, 3, 4, 5, 6, 7, 8, 9," which are used predominantly in the Arab East (\emph{Mashriq}).
    \item The second, the \emph{Ghubar numerals} (known as the Arabic numerals), which were arranged based on angles, comprising: "1, 2, 3, 4, 5, 6, 7, 8, 9," and whose use became prevalent in the Arab West (\emph{Maghrib}).
\end{enumerate}

Crucially, Al-Khwarizmi is credited with adding the \emph{zero (0)} to this numerical system, a contribution that revolutionized mathematics.

Following Europe's engagement with Arabic sciences, beginning in al-Andalus, the West discovered that the dust numerals (Arabic numerals) were more suitable for their use than Roman numerals, owing to the difficulty of employing the latter in addition and subtraction, and their impracticality in multiplication and division~\cite{Harbi2004}.

Credit is thus due to Arab astronomers for transmitting Indian arithmetic~\cite{Saliba2005}, particularly to the Arab scholar al-Khwarizmi, who introduced both Indian and Arabic numerals to the West through his work \textit{al-Sindhind}. These numerals were initially named after him as ``Algorisms''~\cite{Harbi2004}.

The \textit{Z\={\i}j} of Ulugh Beg, known as the ``New Sultanic \textit{Z\={\i}j},'' stands among the most significant astronomical works, in which he compiled and corrected the results of previous observations conducted over twelve years~\cite{Abdullah2016}. Its tables were translated into Latin, and Europeans benefited greatly from them~\cite{AlSheikh2014}.

David A. King (b. 1360 AH/1941 CE) has noted that the \textit{z\={\i}j}s of Ulugh Beg and al-K\=ash\={\i} (d. ca. 832 AH/1428 CE)~\cite{AlMumin1992}, as well as the ``New \textit{Z\={\i}j}'' by Ibn al-Sh\=a\d{t}ir (704 AH/1304 CE--777 AH/1375 CE)~\cite{Kennedy1983}---considered the pinnacle of astronomical science in Syria during the Islamic era---were adapted and simplified, while the \textit{z\={\i}j}s of Cassini and Lalande, translated from French into Turkish and Arabic, failed to replace them~\cite{King1977}.

The Western scholar Liberi (fl. 1838 CE) mentions the existence of a manuscript titled \textit{Liber ysagomarum alchorismi in artem astronomican in Magistro A compositus}, which represents one of the Latin translations in Paris of al-Khwarizmi's book on Indian arithmetic. It contains astronomical tables demonstrating that its author, Pierre Alfonso (1227 AH/1813 CE--1270 AH/1854 CE), based his work upon the tables established by al-Khwarizmi~\cite{Saliba2005}.

Furthermore, a Latin translation of the \textbf{Toledo Tables} (\textit{Z\={\i}j al-T\d{u}lay\d{t}ilah}) by the Arab scholar \textbf{Ab\=u Is\d{h}\=aq al-Zarq\=al\={\i}} (d. ca. 490 AH/1096 CE)~\cite{AlMumin1992} was found in Munich, executed by Gerardus Cremonensis (508 AH/1114 CE--583 AH/1187 CE)~\cite{Saliba2005}.

The \textbf{Alfonsine Tables} (\textit{Tabulae Alphonsinae}) served as the most important link between nascent European astronomy and Islamic science. They were compiled in the second half of the thirteenth century, following the model of al-Zarq\=al\={\i}'s Toledo Tables, at the request of Alfonso X, King of Castile.

The Alfonsine Tables were systematized and became dominant in Europe during the fourteenth century, notably in Paris in 1320 CE, and they surpassed the well-known Toledo Tables of the thirteenth century in many respects.

Their first printed edition dates back to 1443 CE, and this was one of the first astronomical books purchased by \textbf{Copernicus} during his studies in Cracow~\cite{Fahd1981}.

Al-Zarq\=al\={\i} was the first to establish Andalusian astronomy in the eleventh century in Cordoba, and Andalusian astronomy flourished through his works. Al-Zarq\=al\={\i}'s greatness stems from his profound mastery in the use of astronomical instruments.

Al-Zarq\=al\={\i}'s works and his tables, known as the Toledo Tables (compiled after 1068 CE), were the foundation of Andalusian astronomy, remaining the primary reference in astronomy for three centuries in Europe, until the emergence of the Alfonsine Tables. They were the main reference in Paris, London, Toulouse, Pisa, and Marseille from the tenth to the thirteenth century CE.

\begin{figure}[H]
    \centering
    \includegraphics[width=0.5\linewidth]{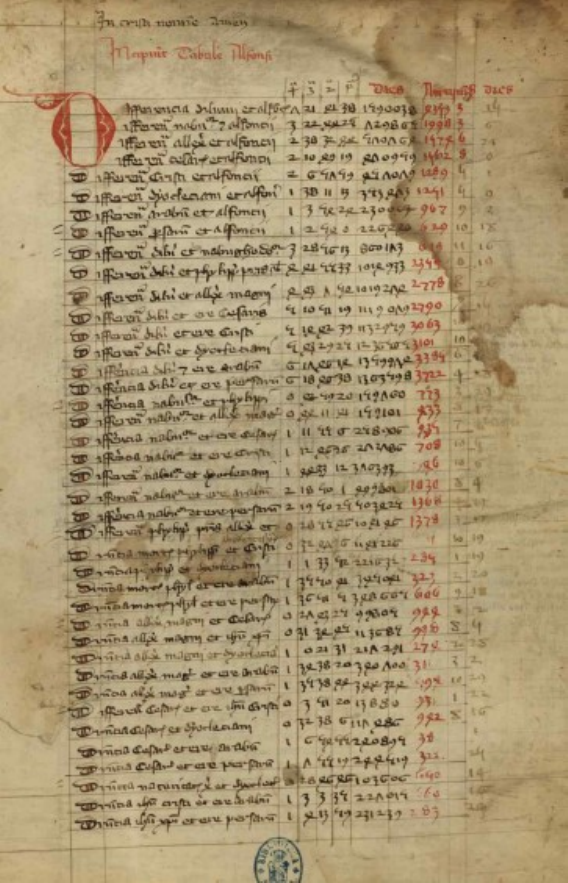}
    \caption{ A snapshot of the tables attributed to Alfonso (fol. 1r).)~\cite{Alfonso1252}.}
    \label{tablets Alphonso}
\end{figure}

In the fourteenth century CE, there were numerous copies of the Toledo Tables, and they remained widespread until the fifteenth century CE, with astronomers continuing to use them alongside the Alfonsine Tables (see Figure~\ref{tablets Alphonso}).

We conclude with the words of Delambre, the author of the history of astronomy:
\begin{quote}
``If you count two or three observers among the Greeks, and then look at the Arabs, you will find a large number of observers among them''~\cite{AlRafei1979}.
\end{quote}
It is therefore not fitting for a scholar to state a fact out of bias for a particular race, as true science and the love of truth are inseparable. Regardless of the assessment of the Greeks' share in the mathematical heritage, the undeniable truth is that they borrowed from the East before the East borrowed from them, and that the people of this East were the ones who entrusted the Europeans with that share, whether large or small, and added to it through refinement and innovation~\cite{AlAkkad2002}.

Regarding the development of astronomy, Arab astronomers played a significant role in the astronomical studies of the Chinese since the thirteenth century, despite the latter's civilizational genius.

The current era of scientific research is closely linked to the study and publication of works related to the history of Arab/Islamic sciences by major Arab scholars, particularly concerning astronomical works. This field has garnered significant attention from Russian historians of science (Mamedbeyli, Rozenfeld, Roganskaya, Bulgakov, Akhmedov, Sergeeva, Matvievskaya, and others) regarding various astronomical and mathematical issues~\cite{AlHamza2022}.
\section{Major Contributions of Orientalists to Arabic Sciences}

Orientalists, particularly \emph{David A. King}, have made significant contributions to the study and understanding of Arabic and Islamic sciences. They worked to bridge the cultural gap between the East and the West through their studies and research. Their efforts contributed to the preservation of the Arab heritage, as well as the documentation and study of Arabic scripts and manuscripts related to Islamic theology (\textit{Kal\=am}) and the Islamic tradition. They also worked on improving the understanding of manuscripts and Arabic texts. Furthermore, their work included the translation and publication of works concerning Islamic schools of thought and sects, which aided in dealing with the diverse general leadership of the world. Their contribution to the development of academic fields related to the Arabic language should also be acknowledged.

David King documented the use of astronomy in the service of Islamic civilization for over a thousand years, through his visits to Arab countries and his related research (one of which was his visit and work at the Institute for the History of Arabic Science in Aleppo, Syria). He published the first comprehensive overview of astronomy in Egypt, Yemen, and Morocco. He also authored a comprehensive book on technological databases in Islamic philosophy, in which he discussed the astronomical precision of prayer times and the determination of the \textit{Qibla} direction. He demonstrated that scholars and philosophers did not rely solely on Greek knowledge but advanced with new and sophisticated endeavors in this field~\cite{King2024}.

These contributions played a role in cultural exchange, significantly influencing our understanding of Arabic and Islamic sciences and helping to preserve and disseminate the Arab heritage to all~\cite{Abdullah2018}.

\section{Results}

\begin{enumerate}
    \item Arab astronomers played a major role in establishing the science of \textit{Z\={\i}j} (astronomical tables) in the West, which was manifested through translation, refinement, and the addition of their own observations.
    \item The author of the Toledo Tables was the Arab scholar al-Zarq\=al\={\i}, not Alfonso, who merely ordered their translation and had them attributed to him, only for them to be later discarded and consequently lost.
\end{enumerate}
\section{Recommendations}

\begin{enumerate}
    \item Highlight the merit of Arab astronomers and their significant role in establishing the science of astronomy in European sources.
    \item The necessity for Arab historians and researchers to collaborate with their Western counterparts to bridge the gap and extract the Arabic sciences from the ruins of the centuries.
\end{enumerate}

\section{Conclusion}

It is clearly evident that Arab and Muslim scholars played a pivotal role in the advancement of Arabic sciences, particularly in the field of astronomy, and continuously contributed to its development in the West. They preserved and developed the heritage of their predecessors, recognizing its value in building modern astronomical knowledge. The role of Orientalists in transmitting this knowledge to the West cannot be overlooked, despite the challenges and criticisms their work faced. The study of the history of Arab and Islamic sciences requires an understanding of the cultural and scientific elements exchanged between the East and the West. Orientalists must appreciate the scientific additions, civilizational roles, and great achievements that the Arab/Islamic civilization has offered to the world.

\section*{Declarations}
\begin{itemize}
	\item Funding: Not Funding.
	\item Conflict of interest/Competing interests: The author declare that there are no conflicts of interest or competing interests related to this study.
	\item Ethics approval and consent to participate: The author contributed equally to this work.
	\item Data availability statement: All data is included within the manuscript.
\end{itemize}

\end{document}